\documentclass[10pt,a4paper,twoside,draft]{amsart}
\usepackage{amssymb}
\usepackage{amsxtra}
\usepackage[mathscr]{eucal}
\usepackage{graphics}
\usepackage{mathtools}
\usepackage{amsrefs}
\usepackage{verbatim}
\usepackage{amscd}

\numberwithin{equation}{section}

\newtheorem{theorem}{Theorem}[section]
\newtheorem{proposition}[theorem]{Proposition}

\newtheorem{corollary}[theorem]{Corollary}

\theoremstyle{definition}

\newtheorem{definition}[theorem]{Definition}
\newtheorem{example}[theorem]{Example}

\newtheorem{problem}[theorem]{Problem}
\newtheorem{remark}[theorem]{Remark}

\newcommand{\al}{\alpha}
\newcommand{\aut}{\operatorname{Aut}}
\newcommand{\be}{\mathbf{e}}
\newcommand{\bj}{\mathbf{j}}
\newcommand{\bk}{\mathbf{k}}
\newcommand{\bm}{\mathbf{m}}
\newcommand{\bn}{\mathbf{n}}
\newcommand{\bs}{\mathbf{s}}
\newcommand{\bu}{\mathbf{u}}
\newcommand{\CC}{\mathbb{C}}

\newcommand{\G}{\Gamma}
\newcommand{\SK}{\mathcal{K}}
\newcommand{\lam}{\lambda}
\newcommand{\lp}{\lam^{(p)}}
\newcommand{\h}{\mathsf{h}}
\newcommand{\m}{\mathsf{m}}
\newcommand{\M}{\mathsf{M}}
\renewcommand{\O}{\Omega}
\newcommand{\On}{\O_n}
\newcommand{\QQ}{\mathbb{Q}}
\newcommand{\Qp}{\QQ_p}
\newcommand{\Pnn}{P_{n\times n}}
\newcommand{\Rd}{\ZZ[u_1^{\pm1},\dots,u_d^{\pm1}]}
\newcommand{\RR}{\mathbb{R}}
\newcommand{\rf}{\rho_f}
\newcommand{\Sig}{\Sigma}
\renewcommand{\SS}{\mathbb{S}}
\newcommand{\TT}{\mathbb{T}}
\newcommand{\U}{\mathsf{U}}
\newcommand{\z}{\zeta}
\newcommand{\ZG}{\ZZ\G}
\newcommand{\ZZ}{\mathbb{Z}}
\newcommand{\zd}{\ZZ^{d}}

\renewcommand{\ge}{\geqslant}
\renewcommand{\le}{\leqslant}
\newcommand{\<}{\langle}
\renewcommand{\>}{\rangle}  

\def\YEAR{\year}\newcount\VOL\VOL=\YEAR\advance\VOL by-1995
\def\firstpage{1}\def\lastpage{1000}
\def\received{}\def\revised{}
\def\communicated{}

\makeatletter
\def\magnification{\afterassignment\m@g\count@}
\def\m@g{\mag=\count@\hsize6.5truein\vsize8.9truein\dimen\footins8truein}
\makeatother

\font\eightrm=cmr8
\font\caps=cmcsc10                    
\font\Caps=cmcsc10 scaled \magstep1   
\font\scaps=cmcsc8

\makeatletter
\@specialpagefalse\headheight=8.5pt
\def\DocMath{{\def\th{\thinspace}\scaps Documenta Math.}}
\renewcommand{\@oddfoot}{\hfill\scaps Documenta Mathematica
    \number\VOL\  (\number\YEAR) \number\firstpage--\lastpage\hfill}
\renewcommand{\@evenfoot}{\ifnum\thepage>\lastpage\hfill\scaps
    Documenta Mathematica \number\VOL\  (\number\YEAR)\hfill\else\@oddfoot\fi}%

\renewcommand{\@evenhead}{%
    \ifnum\thepage>\lastpage\rlap{\thepage}\hfill%
    \else\rlap{\thepage}\slshape\leftmark\hfill{\caps\SAuthor}\hfill\fi}%
\renewcommand{\@oddhead}{%
    \ifnum\thepage=\firstpage{\DocMath\hfill\llap{\thepage}}%
    \else{\slshape\rightmark}\hfill{\caps\STitle}\hfill\llap{\thepage}\fi}%
\makeatother

\def\TSkip{\bigskip}
\newbox\TheTitle{\obeylines\gdef\GetTitle #1
\ShortTitle  #2
\SubTitle    #3
\Author      #4
\ShortAuthor #5
\EndTitle
{\setbox\TheTitle=\vbox{\baselineskip=20pt\let\par=\cr\obeylines%
\halign{\centerline{\Caps##}\cr\noalign{\medskip}\cr#1\cr}}%
	\copy\TheTitle\TSkip\TSkip%
\def\next{#2}\ifx\next\empty\gdef\STitle{#1}\else\gdef\STitle{#2}\fi%
\def\next{#3}\ifx\next\empty%
    \else\setbox\TheTitle=\vbox{\baselineskip=20pt\let\par=\cr\obeylines%
    \halign{\centerline{\caps##} #3\cr}}\copy\TheTitle\TSkip\TSkip\fi%
\centerline{\caps #4}\TSkip\TSkip%
\def\next{#5}\ifx\next\empty\gdef\SAuthor{#4}\else\gdef\SAuthor{#5}\fi%
\ifx\received\empty\relax
    \else\centerline{\eightrm Received: \received}\fi%
\ifx\revised\empty\TSkip%
    \else\centerline{\eightrm Revised: \revised}\TSkip\fi%
\ifx\communicated\empty\relax
    \else\centerline{\eightrm Communicated by \communicated}\fi\TSkip\TSkip%
\catcode'015=5}}\def\Title{\obeylines\GetTitle}
\def\Abstract{\begingroup\narrower
    \parskip=\medskipamount\parindent=0pt{\caps Abstract. }}
\def\EndAbstract{\par\endgroup\TSkip}

\long\def\MSC#1\EndMSC{\def\arg{#1}\ifx\arg\empty\relax\else
     {\par\narrower\noindent%
     2010 Mathematics Subject Classification: #1\par}\fi}

\long\def\KEY#1\EndKEY{\def\arg{#1}\ifx\arg\empty\relax\else
	{\par\narrower\noindent Keywords and Phrases: #1\par}\fi\TSkip}

\newbox\TheAdd\def\Addresses{\vfill\copy\TheAdd\vfill
    \ifodd\number\lastpage\vfill\eject\phantom{.}\vfill\eject\fi}
{\obeylines\gdef\GetAddress #1
\Address #2
\Address #3
\Address #4
\EndAddress
{\def\xs{4.3truecm}\parindent=0pt
\setbox0=\vtop{{\obeylines\hsize=\xs#1\par}}\def\next{#2}
\ifx\next\empty 
     \setbox\TheAdd=\hbox to\hsize{\hfill\copy0\hfill}
\else\setbox1=\vtop{{\obeylines\hsize=\xs#2\par}}\def\next{#3}
\ifx\next\empty 
     \setbox\TheAdd=\hbox to\hsize{\hfill\copy0\hfill\copy1\hfill}
\else\setbox2=\vtop{{\obeylines\hsize=\xs#3\par}}\def\next{#4}
\ifx\next\empty\ 
     \setbox\TheAdd=\vtop{\hbox to\hsize{\hfill\copy0\hfill\copy1\hfill}
                \vskip20pt\hbox to\hsize{\hfill\copy2\hfill}}
\else\setbox3=\vtop{{\obeylines\hsize=\xs#4\par}}
     \setbox\TheAdd=\vtop{\hbox to\hsize{\hfill\copy0\hfill\copy1\hfill}
	        \vskip20pt\hbox to\hsize{\hfill\copy2\hfill\copy3\hfill}}
\fi\fi\fi\catcode'015=5}}\gdef\Address{\obeylines\GetAddress}

\def\LOCAL{\jobname.files}

\begin{document}\frenchspacing

\parskip=\smallskipamount\parindent=0pt

\Title   Mahler's Work and Algebraic Dynamical Systems
\ShortTitle
\SubTitle
\Author  Douglas Lind and Klaus Schmidt
\ShortAuthor
\EndTitle
\Abstract

After Furstenberg had provided a first glimpse of remarkable rigidity
phenomena associated with the joint action of several commuting
automorphisms (or endomorphisms) of a compact abelian group, further key
examples motivated the development of an extensive theory of such
actions.

Two of Mahler's achievements, the recognition of the significance of
Mahler measure of multivariate polynomials in relating the lengths and
heights of products of polynomials in terms of the corresponding
quantities for the constituent factors, and his work on additive
relations in fields, have unexpectedly played important roles in the
study of entropy and higher order mixing for these actions.

This article briefly surveys these connections between Mahler's work and
dynamics. It also sketches some of the dynamical outgrowths of his work
that are very active today, including the investigation of the
Fuglede-Kadison determinant of a convolution operator in a group von
Neumann algebra as a noncommutative generalization of Mahler measure, as
well as diophantine questions related to the growth rates of periodic
points and their relation to entropy.

\EndAbstract
\MSC
Primary: 37A45, 37A15; Secondary 11R06, 11K60.
\EndMSC
\KEY
Algebraic action, entropy, Mahler measure, additive relations in fields.
\EndKEY
\Address Department of Mathematics
University of Washington
Seattle, Washington, 98195
United States
\Address  Mathematics Institute
University of Vienna
Oskar-Morgenstern-Platz 1, A-1090
Vienna, Austria
\Address
\Address
\EndAddress

\section{Dynamical background}\label{sec:dynamical-background}

In order to describe the connections between Mahler's work and dynamical systems,
we have to recall some background information.

Let $X$ be a compact abelian group, and let $\mu$ denote the normalized Haar
measure on $X$, so that $\mu(X)=1$. We write $\aut(X)$ for the group of
continuous algebraic automorphisms of ~$X$. Halmos \cite{Hal} observed 75
years ago that if $A\in\aut(X)$, then the measure $\nu$ defined by
$\nu(E)=\mu(A(E))$ is also a normalized,
translation-invariant measure on $X$, and hence $\nu=\mu$ by uniqueness
of Haar measure. In other words, $A$ preserves the measure~ $\mu$.

\begin{example}
   \label{exam:toral} (\textit{Toral automorphisms})
   Let $\TT=\RR/\ZZ$, and let $\TT^n$ denote the $n$-dimensional torus. Then
   every $A\in GL(n,\ZZ)$ gives an automorphism of $\TT^n$, and all continuous group
   automorphisms of $\TT^n$ arise this way. Hence $\aut(\TT^n)\cong GL(n,\ZZ)$.

   An explicit example to keep in mind is the matrix
   $A=\left[\begin{smallmatrix} 0&1\\1&1\end{smallmatrix}\right]$ acting
   on ~$\TT^2$. Together with its square, the \textit{Arnold cat map} $A^2=\left[\begin{smallmatrix} 1&1\\1&2\end{smallmatrix}\right]$ (cf., e.g., \cite{Arn-Ave}, \cite{Dys}), this toral automorphism has given rise to a vast amount of literature --- mathematical, numerological and phenomenological --- which we cannot explore here, but which certainly makes for fascinating reading.
\end{example}

Rather than exploring intricacies of individual toral automorphisms we shall concentrate here on dynamical systems arising from the simultaneous action of \textit{several} automorphisms of a compact abelian group $X$, and on the somewhat surprising properties of such systems.

Let $\G$ be a countable discrete group (not necessarily abelian). An \textit{algebraic $\G$-action} is a homomorphism $\al\colon\G\to\aut(X)$ for some compact abelian group ~$X$. It is convenient to use exponential notation for $\al$, writing $\al^\gamma$ instead of ~$\al(\gamma)$. In Example \ref{exam:toral}, $\G=\ZZ$ and $\al^k=A^k$.

The interest in algebraic actions of groups other (i.e., bigger) than $\mathbb{Z}$ has its roots in two examples: \textit{Furstenberg's example} \cite{Fur}, consisting of the $\mathbb{N}^2$-action $\alpha $ on $\mathbb{T}=\mathbb{R}/\mathbb{Z}$ generated by the commuting endomorphisms $\times 2$ and $\times 3$, and \textit{Ledrappier's example} \cite{Led}, which will play quite an important role in this article.

\begin{example}
   \label{exam:ledrappier} (\textit{Ledrappier's example}) Consider the
   compact abelian group $Y=(\ZZ/2\ZZ)^{\ZZ^2}$ with coordinate-wise
   addition. Each element $x\in Y$ has the form
   $x=(x_{\bn})_{\bn\in\ZZ^2}$, where each $x_{\bn}\in\ZZ/2\ZZ$, and can
   be thought of as a two-dimensional array of 0's and 1's. There is a
   natural $\ZZ^2$-shift action $\sigma$ on $Y$ defined by
   $(\sigma^{\bm}x)_{\bn}=x_{\bn-\bm}$.   Let $\mathbf{e}_1=(1,0)$ and
   $\mathbf{e}_2=(0,1)$ be the standard basis for ~$\ZZ^2$. Define a
   subgroup $X_L$ of $Y$ by
   \begin{equation}
	\label{eq:ledrappier}
      X_L=\{x\in Y: x_{\bn}+
      x_{\bn+\mathbf{e}_1}+x_{\bn+\mathbf{e}_2}=0
      \text{\enspace for all $\bn\in\ZZ^2$}\}.
   \end{equation}
   This additive condition is clearly shift-invariant, so that we can define
   an algebraic $\ZZ^2$-action $\alpha _L$ on $X$ by restricting $\sigma$ to
   ~$X_L$.
\end{example}

Both these examples are deceptively simple. In Furstenberg's example,
the existence of nonatomic $\alpha $-invariant probability measures $\nu
$ on $\mathbb{T}$ other than Lebesgue measure has remained unresolved
since 1967 and has led to a major new direction of research on
\textit{measure rigidity} of algebraic group actions. In Ledrappier's
example it was the \textit{higher order mixing} properties of the system
which provided the original focus of work by Ledrappier and
others. Another avenue of research opened up with the replacement of the
`alphabet' $\mathbb{Z}/2\mathbb{Z}$ in \eqref{eq:ledrappier} by
$\mathbb{T}$, leading to the notion of a `principal algebraic
$\mathbb{Z}^d$-action' (cf. Example \ref{e:principal}) and, beyond that, to the exploration of \textit{algebraic actions} of arbitrary countable groups ~$\G$.

However, before starting to explore algebraic group actions at that level of generality we return to the more familiar ground of algebraic $\mathbb{Z}^d$-actions with its wealth of examples (see \cite{DSAO} for a detailed account of that theory).

\section{Algebraic $\mathbb{Z}^d$-actions}\label{sec:Zd}

For $\G=\ZZ^d$, the integer group ring $\ZZ\G$ is isomorphic to the ring
$R_d\coloneqq \ZZ[u_1^{\pm1},\linebreak[0]\dots,u_d^{\pm1}]$ of Laurent
polynomials in the commuting variables $u_1,\dots,u_d$. We write a
typical element $f\in R_d$ as
$\sum_{\mathbf{m}\in \mathbb{Z}^d}f_\mathbf{m}\bu^\mathbf{m}$, where
$\bu^{\bm}=u_1^{m_1}\dots u_d^{m_d}$ and $f_\mathbf{m}\in \mathbb{Z}$
with $f_\mathbf{m}=0$ for all but finitely many
$\mathbf{m}\in \mathbb{Z}^d$. When $d=2$, for notational simplicity we
use variables $u$ and $v$ rather than $u_1$ and $u_2$ throughout, so
that $R_2=\ZZ[u^{\pm1},v^{\pm1}]$.

	\begin{example}[Principal and cyclic $\mathbb{Z}^d$-actions]
	\label{e:principal}
There is a natural shift action $\sigma $ of $\mathbb{Z}^d$ on $\TT^{\mathbb{Z}^d}$ given by
$(\sigma^\mathbf{m}x)_\mathbf{n}=x_{\mathbf{n}-\mathbf{m}}$ for every
$x\in \TT^{\mathbb{Z}^d}$ and $\mathbf{m},\mathbf{n}\in \mathbb{Z}^d$.
This definition of the shift map is the opposite of the more traditional
one, but is consistent with how shifts must be defined when the acting
group is noncommutative.
For
$f\in R_d$ define
\begin{equation}
	\label{eq:principal}
   X_f=\Bigl\{x\in\TT^{\mathbb{Z}^d}\colon \sum_{\mathbf{m}\in\mathbb{Z}^d}f_\mathbf{m}
   x_{\mathbf{n}+\mathbf{m}}=0\text{\enspace for all $\mathbf{n}\in \mathbb{Z}^d$}\Bigr\}\subset\TT^{\mathbb{Z}^d}.
\end{equation}
As in Ledrappier's example, this condition on $x$ is invariant under the
shift action $\sigma $ on $\TT^{\mathbb{Z}^d}$, and so the
restriction $\al_f$ of $\sigma $ to $X_f$ gives an algebraic $\mathbb{Z}^d$-action on $X_f$,
called the \textit{principal algebraic $\mathbb{Z}^d$-action defined by $f$}.

For principal actions there is a convenient and very explicit way to
describe the Pontryagin dual $\widehat{X_f}$ and the
$\mathbb{Z}^d$-action dual to ~$\al_f$. The dual group of the cartesian
product $\TT^{\mathbb{Z}^d}$ is the direct sum
$\bigoplus_{\mathbf{m}\in \mathbb{Z}^d}\ZZ$, which as an additive group
is just ~$R_d$. The automorphism dual to the shift-transformation
$\sigma^\mathbf{m}$ is left multiplication by $\bu^\mathbf{m}$ on
~$R_d$. The dual of the subgroup $X_f$ of $\TT^{\mathbb{Z}^d}$ is the
quotient of $R_d$ by the annihilator of $X_f$, which is the principal
ideal ~$fR_d$. Thus $\widehat{X_f}=R_d/fR_d$, which explains the
terminology `principal action'.

If we replace the principal ideal $fR_d$ by an arbitrary ideal $I\subset
R_d$, we obtain the cyclic $R_d$-module $M=R_d/I$ and the corresponding
\textit{cyclic} algebraic $\mathbb{Z}^d$-action $\alpha _{R_d/I}$ on
$X_{R_d/I}=\widehat{R_d/I}$. When $I=f R_d$ is principal, we abbreviate
these to $\al_f$ on $X_f$, as above.
	\end{example}

	\begin{example}[A toral automorphism as a principal $\mathbb{Z}$-action]
	\label{exam:toral-as-principal}
Let $\Gamma =\ZZ$. Then $\ZZ\G$ is isomorphic to $\ZZ[u^{\pm1}]$, the ring of Laurent polynomials in a single variable ~$u$. If $f(u)=u^2-u-1$, it is an instructive little exercise to show that the principal $\mathbb{Z}$-action $(X_f,\al_f)$ is isomorphic to the toral automorphism $A=\left[\begin{smallmatrix} 0&1\\1&1\end{smallmatrix}\right]$ in Example \ref{exam:toral}.

The Arnold cat map $B=A^2=\left[\begin{smallmatrix} 1&1\\1&2\end{smallmatrix}\right]$ in Example \ref{exam:toral} is, of course, also of the form $B=\alpha _M$ for some $R_1$-module ~$M$. Show that this module $M$ is again cyclic.

However, the third power $C=A^3=\left[\begin{smallmatrix} 1&2\\2&3\end{smallmatrix}\right]$ is of the form $C=\alpha _N$ for some $R_1$-module ~$N$ which is \textit{not} cyclic (cf. \cite{DSAO}*{Example 5.3 (2)}). What about $A^n$ with $n>3$? Are any of the corresponding $R_1$-modules cyclic?
	\end{example}

	\begin{example}[Furstenberg's example, revisited]
	\label{e:furstenberg}
Put $\Gamma =\mathbb{Z}^2$ and so  $\mathbb{Z}\Gamma =\ZZ[u^{\pm1},v^{\pm1}]$.
Let $I=\<u-2,v-3\>=(u-2)R_2+(v-3)R_2\subset R_2$ be the nonprincipal ideal generated by $u-2$ and ~$v-3$. As in Example \ref{e:principal} we see that the cyclic $\mathbb{Z}^2$-action $\alpha =\alpha _{R_2/I}$ is the restriction of the shift-action $\sigma $ on $\mathbb{T}^{\mathbb{Z}^2}$ to the closed, shift-invariant subgroup
	\begin{displaymath}
X=X_{R_2/I} =  \bigl\{x\in\TT^{\mathbb{Z}^2}: x_{\mathbf{n}+\mathbf{e}_1}=2x_\mathbf{n}\enspace\textup{and}\enspace  x_{\mathbf{n}+\mathbf{e}_2}=3x_\mathbf{n}\enspace \textup{for every}\enspace \mathbf{n}\in \mathbb{Z}^2\bigr\}.
	\end{displaymath}
If $\pi _\mathbf{0}\colon X\rightarrow \mathbb{T}$ is the projection which sends each $x=(x_\mathbf{n})_{\mathbf{n}\in \mathbb{Z}^2} \in X$ to its zero coordinate $x_\mathbf{0}$, then the diagrams
	\begin{displaymath}
\begin{CD} X @>\alpha ^{\mathbf{e}_1}>>
X\\ @V\pi _\mathbf{0}VV @VV\pi _\mathbf{0}V\\ \mathbb T @>>\times 2> \mathbb T
\end{CD}\qquad \qquad
\begin{CD} X @>\alpha ^{\mathbf{e}_2}>>
X\\ @V\pi _\mathbf{0}VV @VV\pi _\mathbf{0}V\\ \mathbb T @>>\times 3> \mathbb T
\end{CD}
	\end{displaymath}
commute, so that we obtain Furstenberg's example as a factor of $\alpha _{R_2/I}$.
	\end{example}

	\begin{example}[Ledrappier's example, revisited]
	\label{e:ledrappier}
If we identify the additive group $Y=(\mathbb{Z}/2\mathbb{Z})^{\mathbb{Z}^2}$ in Example \ref{exam:ledrappier} with $\{0,\frac12\}^{\mathbb{Z}^2}\subset \mathbb{T}^{\mathbb{Z}^2} = \widehat{R_2}$, then $\widehat{Y}=R_2/2R_2$, and the group $\widehat{X}$ dual to \eqref{eq:ledrappier} is the cyclic $R_2$-module $R_2/\<2,f\>$, where $f(u,v)=1+u+v\in R_2$, and where $\<2,f\>=2R_2+fR_2$ is the nonprincipal ideal generated by $2$ and ~$f$.

For explicit calculations with Ledrappier's example it will be convenient to rewrite this $R_2$-module by viewing $f$ as an element $\tilde{f}$ of the ring $R_2^{(2)} \coloneqq \mathbb{F}_2[u^{\pm1},v^{\pm1}]$ of Laurent polynomials in $u,v$ with coefficients in the prime field $\mathbb{F}_2=\mathbb{Z}/2\mathbb{Z}$, and by identifying $\widehat{X}$ with the $R_2$-module $R_2^{(2)}/ \tilde{f}R_2^{(2)}$.
	\end{example}

	\begin{example}[Ledrappier's example with continuous alphabet]
	\label{e:ledrappier_cont}
If we replace the alphabet $\mathbb{Z}/2\mathbb{Z}$ in Ledrappier's Example \ref{exam:ledrappier} by $\mathbb{T}$, we obtain the closed, shift-invariant subgroup
	\begin{equation}
	\label{eq:ledrappier_cont}
X' = \{x\in \mathbb{T}^{\mathbb{Z}^2}:x_\mathbf{n}+x_{\mathbf{n}+\mathbf{e}_1}+x_{\mathbf{n}+\mathbf{e}_2}=0\enspace \textup{for every}\enspace \mathbf{n}\in \mathbb{Z}^2\}.
	\end{equation}
It is easy to check that the shift-action of $\mathbb{Z}^2$ on $\widehat{X'}$ coincides with the principal $\mathbb{Z}^2$-action $\alpha_f$ on $X_f$, where $f(u,v)=1+u+v$.
	\end{example}

\section{Mixing properties of Algebraic $\mathbb{Z}^d$-actions}\label{sec:mixing}

A little experimentation with equation \eqref{eq:ledrappier} in Ledrappier's Example \ref{exam:ledrappier} yields that
	\begin{equation}
	\label{eq:led2}
x_\mathbf{n} + x_{\mathbf{n}+2^k\mathbf{e}_1} +
x_{\mathbf{n}+2^k\mathbf{e_2}} \equiv 0 \pmod2,\quad k\ge0,
	\end{equation}
so that the coordinates $x_\mathbf{n}$ and $x_{\mathbf{n}+2^k\mathbf{e}_1}$ together determine the coordinate $x_{\mathbf{n}+2^k\mathbf{e}_2}$ of $x$ for every $x\in X$, $\mathbf{n}\in \mathbb{Z}^d$ and $k\ge0$.

Having identified the dual group of Ledrappier's example with
$R_2^{(2)}/ \tilde{f}R_2^{(2)}$ in Example \ref{e:ledrappier}, equation
\eqref{eq:led2} comes as no surprise: since $(g+h)^2=g^2+h^2$ for all
$g,h\in R_2^{{(2)}}$, it follows that
$\tilde{f}(u,v)^{2^k}=1+u^{2^k}+v^{2^k}$ lies in the annihilator
$\tilde{f} R_2^{(2)}$ of $X$ for every $k\ge1$, precisely the content of \eqref{eq:led2}.

In order to appreciate the significance of \eqref{eq:led2}, we recall
that a measure-preserving action
$T\colon \mathbf{n}\mapsto T^\mathbf{n}$ of $\mathbb{Z}^d$ on
a probability space $(Y,\mathscr{T},\nu )$ is \textit{mixing} if
		\begin{displaymath}
\lim_{\|\bn\|\to\infty }\nu (B_1\cap T^{\mathbf{n}}B_2)=\nu (B_1)\nu (B_2)
		\end{displaymath}
for all sets $B_1,B_2\in\mathscr{T}$, where $\|\bn\|$ denotes the euclidean norm of $\bn$. More generally, the action $T$ is \textit{$r$-mixing} with $r\ge 2$ if, for all $B_1,\dots ,B_r\in\mathscr{T}$,
	\begin{equation}
	\label{eq:mixing}
\nu \biggl(\bigcap_{i=1}^rT^{\mathbf{n}_i}B_i\biggr) \longrightarrow \prod _{i=1}^r\nu (B_i)\enspace \textup{as} \enspace \|\mathbf{n}_i-\mathbf{n}_j\|\to\infty \enspace \textup{for}\enspace 1\le i < j \le r.
		\end{equation}

For single measure-preserving automorphisms of probability spaces, the
question of whether mixing implies mixing of every order has been open
for well over 50 years  \cite[p.\ 99]{HalBook}. However, Ledrappier's
example shows that the answer is negative for $\zd$-actions with
$d\ge2$, which we now explain.

It is relatively simple to show that there are no long-range
correlations between pairs of coordinates, so that Ledrappier's example
$(X,\al)$ is mixing. However, \eqref{eq:led2} shows that if $B=\{x\in
X\colon x_{\mathbf{0}}=0\}$, then for all $k\ge1$ we have that
\begin{displaymath}
   B\cap \al^{2^k\be_1}(B)\cap\al^{2^k\be_2}(B)=B\cap\al^{2^k\be_1}(B).
\end{displaymath}
Hence
\begin{equation}
   \label{eqn:led-non-mix}
   \mu\bigl( B \cap \al^{2^k\be_1}(B)\cap\al^{2^k\be_2}(B)\bigr)=
   \mu(B\cap\al^{2^k\be_1}(B)\bigr)=\frac14\ne\frac18=\mu(B)^3
\end{equation}
for all $k\ge1$. Thus Ledrappier's example is not 3-mixing.

In order to reflect the particularly regular way in which higher-order mixing breaks down in Ledrappier's example we introduce a definition.
	\begin{definition}
	\label{d:mixing-set}
Let $T\colon \mathbf{n}\mapsto T^\mathbf{n}$ be a measure-preserving $\mathbb{Z}^d$-action on a probability space $(Y,\mathscr{T},\nu )$. A nonempty finite set $F\subset \mathbb{Z}^d$ is \textit{mixing} if
	\begin{displaymath}
\lim_{k\to\infty }\nu \Bigl(\bigcap_{\mathbf{n}\in F}T^{k\mathbf{n}}(B_{\mathbf{n}})\Bigr)= \prod_{\mathbf{n}\in F} \nu (B_{\mathbf{n}})
	\end{displaymath}
for every collection of Borel sets $B_\mathbf{n}\in\mathscr{T}$, $\mathbf{n}\in F$.
A nonempty finite set $F\subset \mathbb{Z}^d$ is called \textit{nonmixing} if it is not mixing.
	\end{definition}

According to \eqref{eqn:led-non-mix}, Ledrappier's example has the nonmixing set $F=\{\mathbf{0},\be_1,\be_2\}$ of size ~$3$.

More generally, suppose $p\ge2$ is a rational prime, $\mathbb{F}_p=\mathbb{Z}/p\mathbb{Z}$ the corresponding prime field, and $R_d^{(p)}=\mathbb{F}_p[u_1^{\pm1},\dots ,u_d^{\pm1}]$ is the ring of Laurent polynomials in $u_1,\dots ,u_d$ with coefficients in $\mathbb{F}_p$. We write a typical element $\tilde{f}\in R_d^{(p)}$ as $\tilde{f} = \sum_{\mathbf{n}\in \mathbb{Z}^d}\tilde{f}_\mathbf{n}u^\mathbf{n}$ with $\tilde{f}_\mathbf{n}\in \mathbb{F}_p$ for every $\mathbf{n}\in \mathbb{Z}^d$ and denote by $\mathcal{S}(\tilde{f})=\{\mathbf{n}\in \mathbb{Z}^d:\tilde{f}_\mathbf{n}\ne0\}$ the \textit{support} of ~$\tilde{f}$. Exactly as in the brief discussion of Ledrappier's example at the beginning of this section one obtains the following result.

	\begin{proposition}
	\label{p:ledrappier}
Let $I\subset R_d^{(p)}$ be an ideal and let $\alpha _{R_d^{(p)}/I}$ be the cyclic $R_d$-action on the group $X_{R_d^{(p)}/I}$ defined by the ring $R_d^{(p)}/I$, viewed as an $R_d$-module. For every $\tilde{f}\in I$, the set $\mathcal{S}(\tilde{f})$ is nonmixing for $\alpha _{R_d^{(p)}/I}$.
	\end{proposition}

Perhaps surprisingly, such an action $\alpha _{R_d^{(p)}/I}$ may have nonmixing sets which do \textit{not} originate from elements of the ideal ~$I$. For example, if $\tilde{g}(u,v)=1+u+u^2+uv+v^2 \in R_2^{(2)}$, then $\alpha _{R_2^{(2)}/\tilde{g}R_2^{(2)}}$ also has the nonmixing set $F=\{\mathbf{0},\be_1,\be_2\}$ appearing in Ledrappier's example, although the ideal $\tilde{g}R_2^{(2)}$ does not contain any element whose support has cardinality~ $3$. For explanation and details we refer to \cite{KitSch2} and \cite{DSAO}*{Section 28}.

The question of existence -- or nonexistence -- of nonmixing sets for general algebraic $\mathbb{Z}^d$-actions turns out to be intimately connected with a result by Kurt Mahler in his paper \cite{Mahler-additive} on Taylor-coefficients of rational functions. The results in \cite{KitSch} and \cite{Sch-mixing} -- \cite{Sch-variety} show that a mixing algebraic $\mathbb{Z}^d$-action $\alpha $ on a compact abelian group $X$ has nonmixing sets if and only if the dual $R_d$-module $\widehat{X}$ has an associated prime ideal $I\subset R_d$ with the following properties: if $\mathbb{F}=\textup{Quot}(R_d/I)$ is the field of fractions of the integral domain $R_d/I$, and if $G\subset \mathbb{F}$ is the multiplicative subgroup generated by the images in $R_d/I$ of the monomials $u_1,\dots ,u_d \in R_d$, then there exist finitely many elements $a_1,\dots ,a_r$ in $G$ and a nonzero element $(c_1,\dots ,c_r)\in \mathbb{F}^r$, such that
	\begin{equation}
	\label{eq:mahler1}
c_1^{}a_1^m+\dots + c_r^{}a_r^m=1\enspace \textup{for infinitely many}\enspace m\ge1.
	\end{equation}
If the group $X$ is connected, the field $\mathbb{F}=\textup{Quot}(R_d/I)$ in \eqref{eq:mahler1} has characteristic zero for every prime ideal $I$ associated with $\widehat{X}$, and if $X$ is totally disconnected, $\mathbb{F}=\textup{Quot}(R_d/I)$ always has positive characteristic. For Ledrappier's example, the prime ideal in question is $I=\langle 2,1+u+v\rangle $, and the field $\mathbb{F}=\textup{Quot}(R_d/I)$ has characteristic $2$.

In the former case, when $\mathbb{F}$ has characteristic zero, a beautiful $p$-adic argument by Kurt Mahler in \cite{Mahler-additive}*{p. 57} shows that \eqref{eq:mahler1} implies the existence of integers $1\le k<l\le r$ and $b>0$ such that
	\begin{equation}
	\label{eq:mahler2}
a_k^b=a_l^b.
	\end{equation}
By translating this back into our dynamical setting one obtains a contradiction to the hypothesis that $\alpha $ is mixing. This leads to the following conclusion.

	\begin{corollary}[\cite{Sch-mixing}, \cite{DSAO}*{p. 268}]
	\label{c:mixing}
Let $\alpha $ be a mixing algebraic $\mathbb{Z}^d$-action on a compact connected abelian group ~$X$. Then every nonempty finite subset $F\subset \mathbb{Z}^d$ is mixing for ~$\alpha $.
	\end{corollary}

In the latter case, when $\mathbb{F}$ has positive characteristic, David
Masser \cite{KitSch2} proved that \eqref{eq:mahler1} implies a more
complicated relationship between the $a_1,\dots ,a_r$: if $\mathbb{F}$ has
characteristic $p>1$, and if $\overline{\mathbb{F}}_p$ is the algebraic
closure of the prime field $\mathbb{F}_p$, then there exist elements
$b_1,\dots ,b_r$ in the algebraic closure $\overline{\mathbb{F}}$ of
$\mathbb{F}$
and $k,l\ge1$
such that $a^{}_i=b_i^k$ for $i=1,\dots ,r$ and $\{b_1^l,\dots ,b_r^l\}$ is linearly dependent over~ $\overline{\mathbb{F}}_p$.

This result allows one in principle to determine the nonmixing sets of algebraic $\mathbb{Z}^d$-actions on zero-dimensional compact abelian groups.

The story of the connection between mixing properties of algebraic $\mathbb{Z}^d$-actions and additive relations in fields, which turns out to have begun with Mahler's Theorem \cite{Mahler-additive}*{p. 57}, doesn't end with Corollary \cite{Sch-mixing}*{Corollary 2.3}. Subsequent developments, based on remarkable work on $S$-unit equations both in characteristic zero (cf., e.g., \cite{Por-Schl}) and in positive characteristic (cf., e.g., \cite{KitSch2}, \cite{Mas}, \cite{Derksen-Masser1}, \cite{Derksen-Masser2}), clarified the connection between nonmixing sets and the order of mixing. The following result gives a brief and incomplete summary of these later developments.

	\begin{theorem}
	\label{t:converse}
Let $\alpha $ be a mixing algebraic $\mathbb{Z}^d$-action on a compact abelian group ~$X$.
	\begin{enumerate}
	\item
If $X$ is connected, $\alpha $ is mixing of every order \textup{(}\cite{Sch-Ward}*{Corollary 3.3}\textup{)}.
	\item
If $X$ is totally disconnected, then $\alpha $ has nonmixing sets if and only if it does not have completely positive entropy. Furthermore, if $r\ge2$, then $\alpha $ is $r$-mixing if and only if every subset $F\subset \mathbb{Z}^d$ of cardinality $r$ is mixing for $\alpha $ \textup{(}\cite{Mas}*{p. 190}\textup{)}.
	\end{enumerate}
	\end{theorem}

\section{Entropy and Mahler measure}\label{sec:entropy-and-mahler}

Entropy is a numerical invariant of dynamical systems which can be defined for measure-preserving as well as continuous actions. Here our focus will on `topological' entropy, which provides a rough measure of the distortion of the topology of a space under a group action by homeomorphisms of that space.

The exact calculation of entropy, or even a numerical approximation,
is in general difficult. However, computing entropy for algebraic actions
is easier for a very important reason: the homogeneity of an algebraic action means that the calculation of the amount of `distortion' of the space can
be reduced to measuring the distortion of small neighbourhoods of the identity of the group under the action.

In this section we consider an algebraic $\ZZ$-action $\al$ on a compact
abelian group $X$ equipped with Haar measure ~$\mu$.

Let $\mathcal{U}$ be an open neighbourhood of the identity $0_X$
of ~$X$. The set of points in $X$ that remain within $\mathcal{U}$ for the
first $n$ iterates of $\al$ is $\bigcap_{j=0}^{n-1}\al^{-j}\mathcal{U}$, and the
rate of decay of the measure of this set measures how $\mathcal{U}$ changes under the first few elements $\alpha ^1,\dots ,\alpha ^{n-1}$ of the action ~$\alpha $. In order to obtain a scale invariant quantity we consider decreasing sequences of neighbourhoods of the identity and define the \textit{entropy $\h(\al)$ of $\al$} as
\begin{equation}
   \label{eqn:entropy-def}
   \h(\al):=\lim_{\mathcal{U}\searrow \{0_X\}}\,\limsup_{n\to\infty} \,\,-\frac{1}{n}\log
   \mu\biggl(\bigcap_{j=0}^{n-1} \al^{-j}\mathcal{U} \biggr).
\end{equation}

What is not apparent from this definition is the crucial property that
entropy is invariant under measure-preserving conjugacy: if $(X,\al)$ and $(Y,\beta )$ are algebraic $\ZZ$-actions, and if
$\phi\colon X\to Y$ is an invertible measurable map which preserves Haar measure and is \textit{equivariant} in the sense that $\phi\circ\al ^n= \beta ^n\circ\phi$, then $\h(\al)=\h(\beta )$. The
proof of this invariance involves, in particular, establishing the
equality of topological and Haar measure-theoretic entropy for algebraic
actions (see \cite{Wal} for details).

\begin{example}
   \label{exam:k-shift} (\textit{$k$-shift})
   Let $X_k=(\ZZ/k\ZZ)^{\ZZ}$, and $\sigma_k$ be the shift on $X_k$,
   which is called the \textit{$k$-shift}. To
   compute entropy it is enough to consider neighborhoods $\mathcal{U}_r=\{x\in
   X_k\colon x_j=0 \text{\enspace for $-r\le j\le r$}\}$ for large
   ~$r$. Since Haar measure here is product measure,
   \begin{displaymath}
      \mu\biggl(\bigcap_{j=0}^{n-1}\sigma_k^{-j}\mathcal{U}_r^{}\biggr)=\mu\bigl(\{x\in
      X_k \colon x_j=0 \enspace\text{for $-r\le j\le n-1+r$}\}\bigr) =
      \Bigl(\frac{1}{k}\Bigr) ^{n+2r},
   \end{displaymath}
   and hence
   \begin{displaymath}
      -\frac{1}{n}\log
      \mu\biggl(\bigcap_{j=0}^{n-1}\sigma_k^{-j}\mathcal{U}_r^{} \biggr)\to\log k
      \text{\enspace as $n\to\infty$}
   \end{displaymath}
   for every $r\ge1$, and so $\h(\sigma_k)=\log k$.
\end{example}

\begin{example}
   \label{exam:toral-entropy} (\textit{Toral automorphism}) Let
   $A\in GL(r,\ZZ) = \aut(\mathbb{T}^r)$. Then the
   eigenvalues $\lam_1,\dots,\lam_r$ of $A$, listed with multiplicity,
   are all nonzero. The eigenvalues $\lam_i$ for which $|\lam_i|>1$
   control the volume decrease in \eqref{eqn:entropy-def}. If $\mathcal{U}$ is a `nice' small neighbourhood of $0$ in $\mathbb{T}^r$ (like $\prod_{i=1}^r(-\varepsilon ,\varepsilon )$ for small $\varepsilon >0$) then, knowing the Jordan
   form of $A$, it is
   relatively easy to show that we can find positive constants $c_1$ and $c_2$ such that
   \begin{displaymath}
      c_1 n^{-r}\biggl( \prod_{|\lam_i|>1}|\lam_i|^{-n}\biggr)
      \le \mu\biggl(\bigcap_{j=0}^{n-1}A^{-j} \mathcal{U}  \biggr)\le
       c_2 n^{r}\biggl( \prod_{|\lam_i|>1}|\lam_i|^{-n}\biggr).
    \end{displaymath}
    It then follows from the definition \eqref{eqn:entropy-def} that
    \begin{equation}
       \label{eqn:toral-entropy}
       \h(A)=\sum_{|\lam_i|>1}\log|\lam_i|.
    \end{equation}
\end{example}
This formula for entropy was established by Sinai in 1959, shortly after
the introduction of entropy.

\medskip For a polynomial $f(u)\in\CC[u]$, Mahler \cite{Mahler-1-var} defined its
`measure' to be
\begin{displaymath}
   \M(f)=\exp\biggl[\int_0^1 \log|f(e^{2\pi i t})|\,dt \biggr].
\end{displaymath}
It is convenient to introduce the \textit{logarithmic Mahler measure} of
$f$ to be
\begin{displaymath}
   \m(f)=\log \M(f)=\int_0^1 \log|f(e^{2\pi i t})|\,dt.
\end{displaymath}
If we write $f(u)=s\prod_{j=1}^r(u-\lam_j)$, then, as Mahler
observed, Jensen's formula yields that
\begin{displaymath}
   \m(f)=\log s + \sum_{|\lam_j|>1} \log |\lam_j|.
\end{displaymath}

For $A\in GL(r,\ZZ)$ as in Example \ref{exam:toral-entropy}, let
$\chi_A(u)=\det[uI-A]$ be its (monic) characteristic polynomial. The calculation
of entropy in this example can then be expressed as $\h(A)=\m(\chi_A)$,
that is, \textit{entropy equals the logarithmic Mahler measure of the related
characteristic polynomial}.

The first author (DL) observed this around 1970, but viewed it as merely a
curiosity. It does however provide a link between dynamics and the
famous (and still open) Lehmer Problem: does
$\inf\{\m(f)\colon f\in\ZZ[u] \text{\enspace and $\m(f)>0$}\}=0$? The
article \cite{Boyd} by David Boyd in this \textit{Selecta} discusses Lehmer's Problem
in detail.
As shown in \cite{Lind-infinite-torus}, this question is equivalent to
asking whether there are toral automorphisms of arbitrarily small
positive entropy, and also equivalent to asking whether there is an ergodic
automorphism of $\TT^{\ZZ}$ with finite entropy.

The next example shows that $p$-adic fields arise naturally in the study
of algebraic actions, leading to a $p$-adic version of Mahler measure
that is used in calculating entropy for certain actions.

For each rational prime $p$ recall that $\Qp$ denotes the completion of
$\QQ$ with respect to the $p$-adic valuation $|\cdot|_p$, normalized so
that $|p|_p=p^{-1}$. We use the convention that $|\cdot|_\infty$ is the
usual absolute value on $\QQ$, so that $\QQ_\infty=\RR$. Each $\Qp$ is a
locally compact field, and has a normalized Haar measure ~$\mu_p$.

\begin{example}
   \label{exam:2u-3}
   Let $f(u)=2u-3$, and consider the principal algebraic $\ZZ$-action $\al_f$ on
   $X_f$ defined in \eqref{eq:principal}. It is easy to check that $\widehat{X_f}\cong \ZZ[1/6]$ and
   that $\al_f$ is dual to multiplication by $3/2$ on $\ZZ[1/6]$. Then
   locally $X_f$ is $\QQ_2\times\QQ_3\times\RR$, and that in
   this local view $\al_f$ is the diagonal matrix
   \begin{displaymath}
      \begin{bmatrix}
         3/2&0&0\\0&3/2&0\\0&0&3/2
      \end{bmatrix}.
   \end{displaymath}
\end{example}
However, here $3/2$ has different sizes in each component:
\begin{displaymath}
   \Bigl|\frac32\Bigr|_2=2, \quad \Bigl|\frac32\Bigr|_3=\frac13, \quad
   \Bigl|\frac32\Bigr|_\infty=\frac32 .
\end{displaymath}
Reasoning as in Example \ref{exam:toral-entropy}, only those
`eigenvalues' with size greater that 1 contribute to entropy, and so
\begin{displaymath}
   \h(\al_f)=\log 2 + \log \frac32 = \log 3.
\end{displaymath}
Observe also that $\m(f)=\log 2 + \log \frac32$, so that
\begin{displaymath}
   \h(\al_f)=\m(f).
\end{displaymath}
This example combines geometric expansion in $\QQ_\infty$ and arithmetic
expansion in $\QQ_2$ to calculate entropy.

In order to describe a general setting for such results, we introduce
the full solenoid group $\Sig=\widehat{\QQ_d}$, where $\QQ_d$ denotes
the rationals with the discrete topology. If $A\in GL(r,\QQ)$, then $A$
acts via duality on the compact abelian group~ $\Sig^r$. After a series
of papers by several authors dealing with special cases, Yuzvinskii
\cite{Yuz} gave a general formula for the entropy of ~$A$.

\begin{theorem}[Yuzvinskii] Let $A\in GL(r,\QQ)$ have complex eigenvalues
   $\lam_1,\dots,\lam_r$ listed with multiplicity, and let $s$ be the
   smallest positive integer such that $s\chi_A(u)\in\ZZ[u]$. Then
   \begin{displaymath}
      \h(A;\Sig^r)=\log s + \sum_{|\lam_j|>1}\log |\lam_j| .
   \end{displaymath}
\end{theorem}

Yuzvinskii's proof relied heavily on complicated algebra. The role of the
$p$-adics, and the resulting conceptual simplification, was spelled out
in \cite{LindWard}, which we now briefly describe.

Let $A\in GL(n,\QQ)$, and $\overline{\QQ}_p$ denote the algebraic closure
of ~$\Qp$. Then $\chi_A(u)$ factors in $\overline{\QQ}_p[u]$ as
$\chi_A(u)=(u-\lp_1)\dots (u-\lp_r)$. We define the \textit{$p$-adic
logarithmic Mahler measure of $\chi_A(u)$} to be
\begin{displaymath}
   \m_p(\chi_A)=\sum_{|\lp_j|_p>1}\log |\lp_j|_p,
\end{displaymath}
where $|\cdot|_p$ is the (unique) extension of the $p$-adic absolute
value from $\QQ$ to ~$\overline{\QQ}_p$.

We will express the `global' entropy of $A$ acting on $\Sig^r$ as the
sum of `local' entropies, one for each $p\le\infty$. Roughly speaking,
$\Sig^r$ is locally the product $\prod_{p\le\infty}\Qp^r$, and each
factor is preserved by ~$A$. Since entropy adds over products, it follows
that $\h(A;\Sig^r)=\sum_{p\le\infty}h(A;\QQ_p^r)$, where each summand is
the Bowen entropy of a linear map. It is shown in \cite{LindWard} that
$\h(A;\Qp^r)=\m_p(\chi_A)$ for each $p\le\infty$ and that this quantity
vanishes for all but finitely many ~$p$. Hence
\begin{displaymath}
   \h(A;\Sig^r)=\sum_{p\le\infty}\m_p(\chi_A) .
\end{displaymath}
Thus the somewhat mysterious term $\log s$ in Yuzvinskii's formula is simply
the sum of the $p$-adic entropies over $p<\infty$, while the remaining
term is the local entropy at ~$p=\infty$.

\begin{example}
	\label{e:A1}
   Let $A=[3/2]\in GL(1,\QQ)$. Then $s=2$ and $\lam_1=3/2$, and so
   $\h(A;\Sig)=\log3$.
\end{example}

\begin{example}
	\label{e:A2}
   Let $B = \left[\begin{smallmatrix} 0 & -1\\1 &
         6/5 \end{smallmatrix}\right] \in GL(2,\QQ)$. The complex
   eigenvalues of $B$ both have absolute value 1, and so there is no
   geometric contribution to entropy. The only nonzero contribution
   happens in $\QQ_5$, and $\h(B,\Sig^2)=\log5$, providing an interesting
   example where the only expansion is arithmetic.
	\end{example}

\begin{remark}
   Consider the polynomials $f(u)=2u-3$, $g(u)=5u^2-6u+5$ in
   $\mathbb{Z}[u^{\pm1}]$ and define the principal $\mathbb{Z}$-actions
   $\alpha _{f}$ on $X_{f}$ and $\al_g$ on $X_g$ as in
   \eqref{eq:principal}. Then the same calculations as in  Example
   \ref{exam:2u-3} show that $\h(\al_f)=\m(f)=\h(A,\Sig)$ and
   $\h(\al_g)=\m(g)=\h(B,\Sig^2)$. Here the principal algebraic actions
   $(X_f,\al_f)$ and $(X_g,\al_g)$ are equal entropy factors of the
   actions $(\Sig,A)$ and $(\Sig^2,B)$, respectively, appearing in
   Examples \ref{e:A1} and \ref{e:A2}.

   We note that the group $X_f$ coincides with the group $X_{R_2/I}$
   from Furstenberg's Example  \ref{e:furstenberg}, although the
   $\ZZ^2$-action is quite different.
\end{remark}

\section{Algebraic $\zd$-actions and Mahler measure}\label{sec:zd-actions}

In this section we discuss entropy for algebraic $\zd$-actions and the
discovery of its connection with the Mahler measure of polynomials in
several variables.

Let $\al$ be an algebraic $\zd$-action on $X$ with Haar measure
~$\mu$. To define entropy, we simply replace the iterates
$\{0,1,\dots,n-1\}$ used for $\ZZ$-actions in \eqref{eqn:entropy-def}
with the $n$-cube $F_n=\{0,1,\dots,n-1\}^d$ and set
\begin{displaymath}
     \h(\al):=\lim_{\mathcal{U}\searrow \{0_X\}}\,\limsup_{n\to\infty} \,\,-\frac{1}{|F_n|}\log
   \mu\biggl(\bigcap_{\bj\in F_n} \al^{-\bj}\mathcal{U}\biggr).
\end{displaymath}

The crucial property of the $F_n$ is that their boundaries are small
compared with their volumes. More precisely, the obey the \textit{F{\o}lner
condition} that for every $\bk\in\zd$,
\begin{equation}
   \label{eqn:folner}
   \frac{|(F_n+\bk)\triangle F_n|}{|F_n|} \to 0 \text{\quad as $n\to\infty$},
\end{equation}
where $|\cdot|$ denotes cardinality and $\triangle$ denotes symmetric
difference.

\begin{example}
   \label{exam:ledrappier-entropy}
   Recall Ledrappier's example $(X,\mu)$ from Example
   \ref{exam:ledrappier}. Let $F_n=\{0,1,\dots,n-1\}^2$ and
   $\mathcal{U}=\{x\in X\colon x_{(0,0)}=0\}$. Since $\mu$ is
   shift-invariant,
   $\mu\bigl(\bigcap_{\bj\in F_n}\al^{-\bj}\mathcal{U}\bigr)=
   \mu\bigl(\bigcap_{\bj\in F_n}\al^{\bj}\,\mathcal{U}\bigr)$.

   Consider the map $\phi_n\colon X\to(\ZZ/2\ZZ)^{F_n}$
   given by the restriction $\phi_n(x)=x|_{F_n}$. Its image $\phi_n(X)$
   is a subgroup of $(\ZZ/2\ZZ)^{F_n}$ and its kernel is
   $\bigcap_{\bj\in F_n}\al^{\bj}\,\mathcal{U}$. Hence
   \begin{displaymath}
       \mu\biggl(\bigcap_{\bj\in F_n} \al^{\bj}\mathcal{U} \biggr) = \frac{1}{|\phi_n(X)|}.
    \end{displaymath}

    Next we observe that the defining relation
    $x_{(i,j)}+x_{(i,j+1)}+x_{(i+1,j)}=0$ shows that the coordinates of
    $x$ in $F_n$ determine its coordinates in
    $\{(0,0),(1,0),\dots,(2n-1,0)\}$ and conversely, and that the
    coordinates in the latter range may be chosen freely. Hence
    $|\phi_n(X)|=2^{2n}$. Thus as $n\to\infty$ we have that
    \begin{displaymath}
       -\frac{1}{n^2}\log
       \mu\biggl(\bigcap_{\bj\in F_n} \al^{-\bj}\mathcal{U}\, \biggr)=
        -\frac{1}{n^2}\log
       \mu\biggl(\bigcap_{\bj\in F_n} \al^{\bj}\mathcal{U} \,\biggr)
       =
       -\frac{1}{n^2}\log 2^{-2n}\to 0 .
    \end{displaymath}
    A similar argument works for arbitrarily small neighborhoods $\mathcal{U}$ of
    $0_X$, showing that $\h(\al)=0$.
\end{example}

In the spring of 1988 the second author (KS) visited the Institute for
Advanced Study at Princeton. Before leaving for Princeton, KS had been
discussing examples of principal $\mathbb{Z}^2$-actions
$(\alpha _f,X_f)$ for $ f(u,v)\in R_2$ with Tom
Ward, who was a PhD student at Warwick at the time, and who had observed
positivity of entropy for some of these examples. In Princeton, KS
started thinking about positivity of entropy for the `continuous'
version \eqref{eq:ledrappier_cont} of Ledrappier's example, but was
unable to resolve the question.

After Princeton, KS visited Seattle and discussed this problem with DL,
who observed that if the state group is
$C_n=\{0,1/n,2/n,\dots,(n-1)/n\}\subset\TT$ (so that $C_2$ gives
Ledrappier's original example), then each of these `finite
approximations' has zero entropy for the same reason as in the preceding
example. Since $C_n \to\TT$ in some sense, the continuous Ledrappier
example is a limit of zero entropy approximations, and so should also
have zero entropy. DL was so convinced by this reasoning that he bet KS
a Japanese dinner that this was correct.

However, the attempt to turn this intuition into something rigorous ran
into serious difficulties. After fruitless efforts, it began to occur to
the authors that the entropy might be positive after all. It was then
that DL remembered that for principal algebraic $\ZZ$-actions $\al_f$,
entropy equals the logarithmic Mahler measure $\m(f)$ of ~$f$. He wrote a
note to KS which concluded: `And here is a really crazy conjecture: for
general $f$ the entropy should be $\log \M(f)$, where $\M(f)$ is the
Mahler measure of ~$f$.' Motivated by this conjecture, the authors were
subsequently able to obtain the equality
$\h(\alpha _f)=\m(f)\approx 0.3230$ for $f(u,v)=1+u+v$, resulting in a
delicious Japanese dinner for KS. The equality $\h(\alpha _f)=\m(f)$
conjectured by DL was subsequently proved in full generality in
\cite{LSW} and provided the crucial step for computing entropy for
general algebraic $\zd$-actions.

\begin{theorem}[\cite{LSW}]
   Let $0\ne f\in\Rd$. Then the entropy of the associated
   principal algebraic $\zd$-action $\al_f$ is given by $\h(\al_f)=\m(f)$.
\end{theorem}

There are several ways to make this result plausible. We will use the
growth rate of periodic points, which leads to some current research and
open problems.

To motivate what follows, first consider the toral automorphism
$A=\left[\begin{smallmatrix} 0&1\\1&1\end{smallmatrix}\right]$ acting on
~$\TT^2$. It is an instructive exercise to show that the toral
automorphism
$(\TT^2, \left[\begin{smallmatrix} 0&1\\1&1\end{smallmatrix}\right])$ is
the principal algebraic $\ZZ$-action $\al_f$, where $f(u)=u^2-u-1$. Let
$P_n(A)=\{x\in\TT^2\colon A^n x=x\}$, the subgroup of points in $\TT^2$
having period $n$ under ~$A$. To compute $|P_n(A)|$, observe that a point
$x\in\TT^2$ is in $P_n(A)$ iff its lift
$\tilde{x}\in [0,1)^2\subset\RR^2$ satisfies
$(A^n-I)\tilde{x}\in\ZZ^2$. Thus $|P_n(A)|$ equals the number of lattice
points in the parallelogram $(A^n-I)\bigl([0,1)^2\bigr)$, and this is
well-known to be $|\det(A^n-I)|$. Let $\On$ denote the set of $n$th
roots of unity, which is a cyclic subgroup of
$\SS :=\{z\in\CC\colon |z|=1\}$. Let $\lam_1$ and $\lam_2$ be the
eigenvalues of~ $A$. Then
\begin{align*}
  |P_n(A)| = |\det(A^n-I)|&=|(\lam_1^n-1)(\lam_2^n -1)| \\
           &= \biggl| \prod_{\z\in\On}(\lam_1-\z)(\lam_2-\z)\biggr|
             =\prod_{\z\in\On}|f(\z)|.
\end{align*}
Thus we can view
\begin{equation}
   \label{eqn:roots-of-unity}
   \frac{1}{n}\log |P_n(A)| = \frac{1}{n}\sum_{\z\in\On}\log|f(\z)|
\end{equation}
as the logarithmic Mahler measure of $f$ over the subgroup $\On$ of
$\SS$, which we will denote by $\m_{\On}(f)$. Notice that the right-hand
side of \eqref{eqn:roots-of-unity} is a Riemann sum approximation to
$\int_0^1\log|f(e^{2 \pi i s})|\,ds=\m(f)$. Since $\log |f|$ is continuous
on $\SS$, we see that $\m_{\On}(f)\to \m(f)=\h(\al_f)$ as $n\to\infty$, so that
the growth rate of periodic points exists as a limit and equals entropy.

The convergence of $\m_{\On}(f)$ to $\m(f)$ is much more delicate if $f$
has roots on $\SS$, for example if $f(u)=u^4+4u^3-2u^2+4u+1$. Then
$\log|f|$ has logarithmic singularities on $\SS$, and the value of $\log|f(\z)|$
for some $\z\in\On$ could be extremely negative should $\z$ be very
close to a root $\lam$ of $f$, or equivalently, if $|\lam^n-1|$ is very
small. However, a deep diophantine result of Gelfond \cite{Gelfond} says that if
$\lam\in\SS$ is an algebraic number, then for every $\varepsilon>0$
there is a constant $C>0$ such that $|\lam^n-1|>Ce^{-\varepsilon n}$,
and using this one can show the convergence $\m_{\On}(f)\to \m(f)$ for
all $f\in\ZZ[u^{\pm1}]$ that have no roots that are roots of unity.

Let us use this approach on the continuous Ledrappier example, i.e., the
principal algebraic $\ZZ^2$-action $\al_f$, where $f(u,v)=1+u+v$. For
simplicity we start with `square' sublattices $n\ZZ\times
n\ZZ=n\ZZ^2\subset\ZZ^2$. Define
\begin{displaymath}
   \Pnn(\al_f):=\{x\in X_f\colon\al ^{n\bj}x = x\text{\quad for all $\bj\in\ZZ^2$}\}
\end{displaymath}
to be the subgroup of all points in $X_f$ fixed by iterates in $n\ZZ^2$.
A calculation similar to the $1$-dimensional case above suggests that
\begin{equation}
   \label{eqn:2-dim}
   |\Pnn(\al_f)|=\prod_{(\xi,\z)\in\On^2}|f(\xi,\z)|,
\end{equation}
so that
\begin{equation}
   \label{eqn:2-dim-RS}
   \frac{1}{n^2} \log|\Pnn(\al_f)|=\frac{1}{n^2}\sum_{(\xi,\z)\in\On^2} \log|f(\xi,\z)|
\end{equation}
is a Riemann sum approximation to $\m(f)$, and we would therefore expect
that $(1/n^2)\log|\Pnn(\al_f)|\to\h(\al_f)$ as $n\to\infty$.

However, there is a serious problem. If $\omega=e^{2\pi i/3}$, then
$f(\omega,\omega^2)=0=f(\omega^2,\omega)$. Thus if $3\mid n$, two of the
summands in the right-hand side of \eqref{eqn:2-dim-RS} are $\log 0
=-\infty$, and the product in \eqref{eqn:2-dim} equals 0. Dynamically,
what is happening is that when $3\mid n$, the subgroup $\Pnn(\al_f)$ is
no longer finite, but rather a finite union of cosets of a
2-dimensional torus. The solution to this situation is to count the
\textit{connected components} of $\Pnn(\al_f)$, which corresponds to
ignoring those points in $\On^2$ where $f$ vanishes. Thus we
\textit{define}
\begin{displaymath}
   \m_{\On^2}(f):= \frac{1}{n^2} \sum_{(\xi,\z)\in\On^2,\,\,f(\xi,\z)\ne0}
   \log|f(\xi,\z)|.
\end{displaymath}
With this convention, one can prove that $\m_{\On^2}(f)\to\h(\al_f)$,
i.e., that the growth rate of periodic components exists as a limit and
is equal to entropy. For more information about $\m(1+u+v)$ see
\cite{Boyd} in this \textit{Selecta}.

For general $f\in\Rd$ a crucial role is played by the \textit{unitary
variety} of $f$, defined as $\U(f):=\{\bs\in\SS^d\colon f(\bs)=0\}$. As
suggested by the discussion above, if $K$ is a finite subgroup of
$\SS^d$, we define
\begin{equation}
   \label{eqn:finite-meas}
   \m_K(f):=\frac{1}{|K|} \sum_{\bs\in K\smallsetminus
   \U(f)}\log|f(\bs)|.
\end{equation}
We use the notation $K\to\infty$ to mean that Haar measures on the finite
subgroups $K$
converge weakly to Haar measure on $\SS^d$.

\begin{problem}
   \label{prob:finite-conv}
   If $0\ne f\in\Rd$, does
   \begin{equation}
      \label{eqn:converge}
      \m_K(f)\to \m(f) \text{\quad as\quad $K\to\infty$}?
   \end{equation}
\end{problem}
In other words, do the Riemann sums for $\log |f|$ over finite subgroups
of $\SS^d$ (modified to avoid values of $-\infty$) converge to
$\int_{\SS^d}\log|f|$?

In \cite{DSAO}  KS had shown that the answer to Problem \ref{prob:finite-conv} is `yes' if $\h(\alpha _f)<\infty $, and if one replaces the limit in \eqref{eqn:converge} by `lim sup'. In \cite{LSV1} it was shown that the answer is `yes' if $\U(f)$ is
finite, and then in \cite{LSV2} that the  answer is also `yes' if the real
dimension $\dim \U(f)$ of $\U(f)$ is less than or equal to ~$d-2$. Both
papers use dynamical ideas, in particular they use homoclinic points for
the action to create sufficient many periodic components.

\begin{example}
   Consider a three variable version of Ledrappier's example with
   continuous alphabet, defined by $f(u,v,w) = 1+u+v+w$. Here the
   unitary variety $\U(f)\subset\SS^3$ is a union of three circles, each given by
   setting one of the variables equal to ~$-1$. Since these circles are
   cosets of 1-dimensional subgroups, the dimension of the connected
   component of the identity for points fixed by all iterates in
   $n\ZZ^3$ grows linearly in~ $n$. A subtlety in the proof of
   convergence in \cite{LSV2}, as illustrated in this example, is that
   although the dimension of the connected components grows, the linear
   constraint here can be used to show that these components do not contribute to
   entropy.
\end{example}

Convergence in \eqref{eqn:converge} is also a diophantine problem,
essentially asking how close points in $K$ can come to ~$U(f)$. Using
ideas involving diophantine analysis, Vesselin Dimitrov \cite{Dim1}
gave a completely different proof of \eqref{eqn:converge} in his 2017
PhD thesis from Yale, again under the assumption that $\dim
\U(f)\le d-2$. Very recent work of Habegger \cite{Hab} on special points
near definable sets, whose proof relies on logic and \textit{O}-minimal sets, was
used by Dimitrov to provide yet a third proof, quite different from the
previous two, again only in the case $\dim \U(f)\le d-2$.

All three of these proofs fail when $\dim\U(f)=d-1$. The following
example illustrates the difficulties.

\begin{example}
   \label{exam:codim-1}
   Let $f(u,v)=3-u-u^{-1}-v-v^{-1}$. Here $\U(f)$ is a 1-dimensional
   oval in $\SS^2$, so that $\log|f|$ has a 1-dimensional set of
   logarithmic singularities. There are only four points on $\U(f)$ both
   of whose coordinates are roots of unity: $(\omega^{\pm1},1)$ and
   $(1,\omega^{\pm1})$, where $\omega=e^{\pi i/3}$ (see \cite[Example 4.4]{LSV2}).
   Whether $\m_K(f)\to\m(f)$ as $K\to\infty$ is still open. But for
   `square' subgroups $K=\On^2$, Dimitrov \cite{Dim2} has very recently shown
   convergence using quite difficult diophantine arguments, and that this
   holds for all nonzero $f\in\Rd$ using `square' subgroups ~$\On^d$.
\end{example}

We can put Problem \ref{prob:finite-conv} into a more general context as
follows. Let $\SK$ denote the set of all compact subgroups of
~$\SS^d$. For $K\in\SK$ let $\mu_K$ be Haar measure on~ $K$. In analogy
with \eqref{eqn:finite-meas}, define the logarithmic Mahler measure of
$0\ne f\in\Rd$ over $K$ to be
\begin{displaymath}
   \m_K(f) := \int_{K\smallsetminus\U(f)}\log|f(\bs)|\,d\mu_K(\bs),
\end{displaymath}
which agrees with our earlier definition when $K$ is finite.

Now $\SK$ is a compact metric space with respect to the Hausdorff metric
on compact subsets of $\SS^d$. Lawton \cite{Law} showed that the
function $K\mapsto\m_K(f)$ is continuous on the closed subset of $\SK$
consisting of all subgroups having dimension at least 1. For example,
the 1-dimensional subgroups $K_n=\{(s,s^n)\colon s\in\SS\}$ of $\SS^2$ converge to
$\SS^2$ in the Hausdorff metric, and so $\m\bigl(f(u,u^n)\bigr)\to
\m\bigl(f(u,v)\bigr)$ as $n\to\infty$, where $f(u,u^n)$ is considered
as a polynomial in $\ZZ[u^{\pm1}]$. Boyd \cite{Boyd} discusses the rate
of this convergence, and in particular how fast $\m(1+u+u^n)$ converges
to $\m(u,v)$.
We can therefore reformulate Problem \ref{prob:finite-conv} in terms of continuity
of Mahler measures on compact subgroups of ~$\SS^d$.

\begin{problem}
   Fix $0\ne f\in\Rd$. Is the function $K\mapsto \m_K(f)$ continuous on $\SK$?
\end{problem}

\section{Algebraic actions of noncommutative
groups}\label{sec:noncommutative}

For a general countable group $\G$ we denote by $\ZZ\G$ the integral
group ring of ~$\G$, where $\G$ is written multiplicatively. A typical
element in $\ZZ\G$ has the form $f = \sum_{\gamma\in\G}f_\gamma \gamma$,
where each $f_\gamma\in\ZZ$, and where all but finitely many of the
$f_\gamma$ vanish. Multiplication in $\mathbb{Z}\Gamma $ is carried out
in the obvious way to extend multiplication in ~$\G$. If $M$ is a
countable left module over $\ZZ\G$, regarded as a discrete abelian group
under addition, then the Pontryagin dual $X_M=\widehat{M}$ is a compact
abelian group, and we obtain an algebraic $\Gamma $-action $\alpha _M$
on $X_M$ by setting, for every $\gamma \in \Gamma $, $\alpha _M^\gamma $
equal to the automorphism of $X_M$ dual to left multiplication by
$\gamma$ on $M$.

Since every algebraic $\Gamma $-action arises in this manner, there is a
1-1 correspondence between left $\mathbb{Z}\Gamma $-modules and
algebraic $\Gamma $-actions, exactly as for $\Gamma =\mathbb{Z}^d$
(except that we now have to be careful about `left' and
`right'). However, if $\Gamma $ is noncommutative, much less is known
about (left) ideals in and (left) modules over $\mathbb{Z}\Gamma $ than
is the case for $\Gamma =\mathbb{Z}^d$, and even for principal algebraic
$\Gamma $-actions (defined by complete analogy with the principal
$\mathbb{Z}^d$-actions Example \ref{e:principal}) our understanding of
their dynamical properties is rudimentary.

In order to fix notation we again write $\sigma $ for the \textit{left} shift-action
$(\sigma^\theta x)_\gamma=x_{\theta^{-1}\gamma}$ of $\Gamma $ on $\mathbb{T}^\Gamma $, and consider, for $f\in\ZZ\G$, the closed, shift-invariant subgroup
\begin{equation}
	\label{eq:principal Gamma}
   X_f=\Bigl\{x\in\TT^\G \colon \sum\nolimits_{\gamma\in\G}f_\gamma
   x_{\theta\gamma}=0\text{\enspace for all $\theta\in\G$}\Bigr\}\subset\TT^\G.
\end{equation}
As in Example \ref{e:principal} we call the restriction $\alpha _f$ of $\sigma $ to $X_f$ the \textit{principal algebraic $\G$-action defined by $f$}. The shift-transformation
$\sigma ^\gamma $ on $\mathbb{T}^\Gamma $ is again dual to left multiplication by $\gamma $ on $\ZZ\G$, and the dual of the subgroup
$X_f\subset \TT^\G$ is the quotient of $\ZZ\G$ by
the left principal ideal $\ZZ\G f$, that is,
$\widehat{X_f}=\ZZ\G/\ZZ\G f$.

One afternoon about fifteen years ago Wolfgang L\"{u}ck presented a copy
of his new book \textit{$L^2$-Invariants: Theory and Applications to
Geometry and $K$-Theory} to his colleague Christopher Deninger at the
University of M\"{u}nster. Not quite knowing what to do with a thick
book on unfamiliar topics, Deninger started flipping through its pages
randomly. By chance, he stumbled on an example showing that the Mahler
measure of $f\in\Rd$ equals the Fuglede-Kadison determinant of an
associated convolution operator $\rf$ on $\ell^2(\zd)$. Knowing the
connection between entropy of algebraic $\zd$-actions and Mahler
measure, Deninger realized that the convolution operator approach, which
is easily generalized to arbitrary countable groups $\G$, might give a
way to compute entropy for principal algebraic $\G$-actions. He was able
to show in \cite{Den} that his idea works for \textit{expansive}
principal actions of at least a restricted class of amenable groups, a
breakthrough that initiated the serious study of algebraic actions of
arbitrary countable discrete groups. His insight can be viewed as a way
to define Mahler measure for polynomials in noncommuting variables.

Before 2010 it was widely believed that entropy theory for group actions
was restricted to actions by amenable groups. This changed when Lewis
Bowen \cite{BowenInv} introduced radically new ideas that ultimately
allowed an extension of entropy theory to a much larger class of groups
called sofic groups, those having a certain kind of finite
approximation. To date there is no known example of a countable group
that is not sofic.

In this section we briefly sketch these two major developments in
dynamics, and how they combined recently in a comprehensive entropy
theory for algebraic actions.

Let $0\neq f\in\Rd$. Define $f^*(u_1,\dots,u_d) := f(u_1^{-1}, \dots,
u_d^{-1})$. Regard $f$ as a function on $\SS^d$, so that
$f^*=\bar{f}$, the complex conjugate of ~$f$.  Consider the
multiplication operator $T_f$ on $L^2(\SS^d)$ given by
$T_f(\phi)=f\cdot\phi$ for $\phi \in L^2(\SS^d)$. Then
\begin{displaymath}
   T_f^* T_f^{}=T_{f^*} T_f^{}=T_{|f|^2}=T_{|f|}^2 ,
\end{displaymath}
and the spectral measure $\mu_{|f|}$ of $T_{|f|}$ is the push-forward of
Lebesgue measure on $\SS^d$ under the map $|f|$, so that $\mu_{|f|}$ is
supported on the real interval $[0,\|f\|_\infty]$. Fuglede and Kadison
\cite{FK} introduced a notion of determinant for certain classes of
operators that include $T_{|f|}$. We then calculate, using their
definition and  change of
variables, that
\begin{displaymath}
   \det T_{|f|}:= \exp \biggl[ \int_0^{\infty} \log t \,d\mu_{|f|}(t)\biggr] =
   \exp \biggl[\int_{\SS^d}\log|f(\bs)|\,d\mu(\bs)\biggr] =\M(f).
\end{displaymath}
This is the fact that Deninger came across in L{\"u}ck's book.

The Fourier transform gives an isomorphism from $L^2(\SS^d)$ to
$\ell^2(\zd)$, and under this isomorphism the multiplication operator
$T_f$ is mapped to the convolution operator $\rf$ on
$\ell^2(\zd)$. Concretely, if we view $w\in\ell^2(\zd)$ as a formal sum
$w=\sum_{\bn\in\zd}w_{\bn}\bu^{\bn}$, then $\rf(w)=w\cdot f$, extending the customary multiplication of polynomials. The
connection with $\al_f$ is provided by the observation that if points
$t\in\TT^{\zd}$ are similarly regarded as formal sums
$t=\sum_{\bn\in\zd}t_{\bn}\bu^{\bn}$, then $X_f=\ker
\rho_{f^*}$. Deninger realized that since
\begin{displaymath}
   \h(\al_f)=\m(f)=\log \det T_{|f|}=\log \det \rf,
\end{displaymath}
the calculation of entropy could be phrased entirely in terms of
convolution operators. This avoids the use of Fourier transforms, and
suggests a general way to deal with principal actions of noncommutative
groups.

Now let $\G$ be a general discrete countable group. As above, consider
points $w\in\ell^2(\G)$ as formal sums $\sum_{\gamma\in\G}
w_{\gamma}\gamma$. For $f\in\ZZ\G$ there is the convolution operator
$\rf$ on $\ell^2(\G)$ given by $\rf(w)=w\cdot f$. The weak operator
closure of the set of complex combinations of these convolutions
operators is called the \textit{group von Neumann algebra}
$\mathcal{L}\G$ of ~$\G$. For $U\in\mathcal{L}\G$, the functional
calculus gives an operator $|U|\in\mathcal{L}\G$ with $U^*U=|U|^2$. The positive
self-adjoint operator $|U|$ has a spectral measure $\mu_{|U|}$
supported on $[0,\|U\|]$. Fuglede and Kadison defined $\det U$ by
\begin{displaymath}
   \det U :=\int_0^\infty\log t \, d\mu_{|U|}(t).
\end{displaymath}
One deep result in this theory is that $\det(UV)=(\det U)(\det V)$. This
definition applies to the operators $\rf$, and we abbreviate $\det \rf$
to $\det f$.

Deninger \cite{Den} showed that his idea worked for amenable groups
having special kinds of F{\o}lner sequences. A series of improvements by
several authors culminated in the definitive result for principal
algebraic actions of amenable groups by Hanfeng Li and Andreas Thom
\cite{LiThom}: Let $\G$ be an amenable group and $f\in\ZG$; if $\rho_f$
is injective on $\ell^2(\G)$, then $\h(\al_f)=\log\det f$, and otherwise
$\h(\al_f)=\infty$. A consequence is that
$\h(\al_{f^*})=\h(\al_{f^{}})$, which is highly nontrivial since there
is no obvious dynamical connection between $\al_{f^{}}$ and $\al_{f^*}$
when $\G$ is noncommutative.

A concrete example of noncommutative $\G$ is the discrete Heisenberg
group $\mathbb{H}$, the group generated by $u$, $v$, and $w$ with
relations $uw=wu$, $vw=wv$, and $vu=wuv$. Even for this simplest
infinite noncommutative group there are many open
problems, e.g., characterize those $f\in\ZZ\mathbb{H}$ for which
$\h(\al_f)=0$, or determine the higher order mixing properties of
principal $\mathbb{H}$-actions. For a comprehensive survey of  what is
currently known about algebraic
$\mathbb{H}$-actions see~ \cite{LS-Heis}.

At roughly the same time, the extension of entropy theory to sofic
groups was undergoing vigorous development, with algebraic actions
providing important and guiding examples. A lucid and systematic account
is contained in the recent book by David Kerr and Hanfeng Li
\cite{KerrLi}. This book describes a profound shift in viewing
$\G$-actions, from the traditional `internal' view using objects
within the space being acted upon to an `external' view using finite
models of the action. In this way F{\o}lner sets and amenability are
avoided, but at the cost of more abstract and complicated machinery,
whose implications are still being worked out.

With the ability to define entropy for principal $\G$-actions for
sofic $\G$, and a viable candidate $\log\det f$ for its value, these two
strands of dynamical progress culminated in the definitive theorem by
Ben Hayes \cite{Hay}, who showed that if $\G$ is sofic and
$f\in\ZZ\G$, then $\h(\al_f)=\log\det f$ provided that $\rf$ is
injective on $\ell^2(\G)$, and is equal to $\infty$ otherwise.

The chain of events set in motion by the discovery in \cite{LSW} that
entropy equals logarithmic Mahler measure for algebraic $\zd$-actions
has led to a remarkable level of generality in the entropy theory of
algebraic actions. However, other dynamical properties of such actions,
like mixing, positivity of entropy, or the Bernoulli property, still remain rather mysterious as soon as one leaves the comfortable world of $\mathbb{Z}^d$-actions.

\Addresses
\end{document}